\documentclass{elsart}
\usepackage{latexsym}
\usepackage{amsfonts}
\usepackage{amsmath}
\usepackage{amssymb}
\def\E{E}

\begin{document}
\begin{frontmatter}
\title{M\"obius Energy of Thick Knots}

\author{Eric J. Rawdon\thanksref{Rawdon}}
\address{Department of Mathematics and Computer Science,\\
Duquesne University,\\ Pittsburgh, PA 15282, USA\\ Email:
rawdon@mathcs.duq.edu}

\author{Jonathan
Simon\thanksref{corresp}\thanksref{Simon}}
\address{Department of Mathematics,\\ University of
Iowa,\\ Iowa City, IA 52242, USA\\ Email:
jsimon@math.uiowa.edu}

\thanks[corresp]{Corresponding author. (J. Simon) Tel:
319-335-0768,  fax: 319-335-0627.}

\thanks[Rawdon]{Research supported by Chatham College,
the University of Iowa, and NSF grant \#DMS0074315}
\thanks[Simon]{Research supported by NSF grant
\#DMS9706789}

\begin{abstract} 
The M\"obius energy of a knot is an
energy functional for smooth curves based on an idea of
self-repelling.  If a knot has a thick tubular 
neighborhood, we would intuitively expect the energy to
be low. In this paper, we give explicit bounds for energy
in terms of the ropelength of the knot, i.e.{} the ratio
of the length of a thickest tube  to its radius.
\end{abstract}

\end{frontmatter}

\section{Introduction}

In this paper, we exhibit a bound for the M\"obius energy
of a knot, in terms of the amount of ``rope'' needed to
make the knot.  This is the energy introduced in 
\cite{O1} and studied extensively in
\cite{FHW,KS}.  We follow an overall approach suggested
by G. Buck for showing that any inverse-square knot
energy should be bounded by the $4/3$ power  of
the ropelength.

We and other researchers have defined a number of
different energy functions for (smooth or polygonal) knots
\cite{O1,FHW,KS,BO1,BS,BS2,BS3,DEJ,Fu,Lo,O3,O2,O4,Si3,Si4,Si2}
based on the idea of inverse-square repelling energy (so
these would correspond to inverse-cube ``forces'').  (See
also \cite{Mo} for a different approach).  Roughly, these
energies are defined in terms of integrals over the curve
$K$
 $$\int\nolimits_{x \in K}\int\nolimits_{y \in
K}{\frac{\Box}{|x-y|^2}}\;\d x\;\d y\;\;.$$ Here  $\Box$
is a placeholder for any of several kinds of terms that
make the integral not give too much weight to the
repelling of points that are close to each other in the
sense of arclength along the knot (so the improper
integral will converge).  
The same purpose is accomplished by defining energies 
as integrals of
differences \cite{O1,FHW},
$$\int\nolimits_{x \in K}\int\nolimits_{y \in
K}{\frac1{|x-y|^2}}\;\;
-\;\;{\frac1{|s-t|^2}}\;\;\;\d x\;\d y\;.$$
Here $s,t$
lie on a line or circle used for a unit-speed parameterization
of the curve $K$ with $s\mapsto x$ and $t\mapsto y$.  
The energy we study in this paper can be defined
either way \cite{O1,FHW,KS} and we shall use the latter.

In addition to viewing a curve as self-repelling, one
also can view it as self-excluding, and define the {\it
ropelength}
 energy $\E_L(K)$ \cite{BO1,BS2,BS3}: this is the ratio of
the arclength of the curve to the maximum radius of a
uniform non self-intersecting tube  along the curve 
\cite{LSDR}, i.e.{} the ratio of length to radius of the
rope (see Section \ref{thelemmas} for precise
definitions).  Variations on thickness are developed in
\cite{Mo,devrthickness,DEJ2,devrthickness2,maddocks,KS2,meideal,mine}.  

We showed in \cite{BS3} that the {\it normal} energy
$\E_N(K)$, which discounts tangential self-repelling, 
and the  {\it symmetric} energy $\E_S(K)$, which models 
self-radiation of a filament \cite{B2}, are bounded by
the ropelength.  These energies, in turn, dominate the
number of crossings in any regular projection of the
knot.  The inequalities are of the form shown below in
(\ref{longineq}).  Here $\mathrm{acn}(K)$ is the {\it
average crossing number}, that is the average, over all
spatial directions, of the number of crossings seen from
each direction.  This, in turn, is larger than the
crossing number
$\mathrm{cr}(K)$, which is the minimum over all regular
projections, and
$\mathrm{cr}[K]$, which considers all $K$ in a given knot
type.

\begin{equation} 4\pi\,\mathrm{cr}[K] \leq 4\pi\,
\mathrm{cr}(K)< 4
\pi\, \mathrm{acn}(K) \leq
\E_S(K) \leq \E_N(K) \leq c\,\E_L(K)^{4/3}
\label{longineq}
\end{equation}

The coefficient $c$ varies with different proofs; the
exponent $4/3$ is sharp.  A related idea is the
{\it writhe} of a knot,
$\mathrm{wr}(K)$,   which is the average over all spatial
directions of the {\it signed} crossing numbers.  Since
$\mathrm{wr}(K) \leq \mathrm{acn}(K)$, we get the same
bound on 
$4\pi\, \mathrm{wr}(K)$.  But using a different analysis
of writhe, in terms of vector fields flowing in tubes
around the knot, it is shown in  \cite{CDG} that
$\mathrm{wr}(K) \leq {\frac14}
\E_L(K)^{4/3}$, which is a lower coefficient than we have
for $\mathrm{acn}(K)$. The
coefficient, approximately $5$, that we obtain in this
paper is lower than in
\cite{BS3}, because we do a more subtle analysis. While
coefficients might be improved, the exponent
$4/3$ in  inequality (\ref{longineq}) is sharp,
 because of examples \cite{B3,CKS} where crossing number
grows like the
$4/3$ power of ropelength. Similarly, because
M\"obius energy also bounds crossing number
\cite{FHW}, the same examples
 show that the  exponent $4/3$ in Theorem
\ref{main} of this paper is sharp.

Our main result  is that the M\"obius energy
 is bounded by  the $4/3$
power of the ropelength (with coefficient
$\approx 5$).   This and the results cited above
support our belief that ropelength is the fundamental
measure of knot complexity: given a bound on ropelength,
one can find a bound on any other given knot invariant.

We state the theorem for the version of the energy that
equals $4$ on a circle, $$\E_{O4}(K) = \int\nolimits_{x \in
K}\int\nolimits_{y \in K} {  {\frac1{|x-y|^2}} - {  
{\frac1{\mathrm{arc}(x,y)^2}}   }}\;\;,$$  where
$\mathrm{arc}(x,y)$ denotes the minimum arclength along
the curve
$K$ between $x$ and $y$.  

We use one kind of analysis
to bound the energy contribution coming from pairs of
points $(x,y)$ where
$x$ and $y$ are near each other in arclength along the
knot, and a different analysis for pairs where $x$ and $y$
are relatively far apart.  The analysis of proximal pairs
will be special for each particular energy function, and
yields typically a {\it linear} bound on energy in terms
of ropelength.  The analysis of distal pairs is identical
for the various energy functions: with little work, we
can get a bound on energy that is {\it quadratic} in
ropelength; and with more work, we obtain a bound that is
$4/3$ power in the ropelength. We add the proximal and
distal contributions to get the overall bound.

In a recent book \cite{SKK}, many ideas of energy and
thickness for knots are discussed.  Many of us interested
in such ideas thought that knots should have ``ideal''
forms.  It has turned out that the knot conformations
that minimize the various ideas of energy differ from one
another.  The  similarities and differences are equally
provocative.

\section{The Lemmas}
\label{thelemmas}

The first lemma is a version of a theorem of Schur.  This
is taken from \cite{Ch}, with the slight adjustment that
we do not need to consider general planar convex curves as
the reference curves, just circles.
 
\begin{lem}  Let $K$ be a $C^2$ smooth curve in $\Rset^3$
whose curvature everywhere is
$\leq$ some number $k$.  Let $C$ be a circle of curvature
$k$, i.e.{} of radius $r = {\frac1k}$.   Let $x,y \in K$ and $s,t
\in C$ 
be any points
such that $\mathrm{arc}(x,y)=\mathrm{arc}(s,t)
\leq \pi r$.  Then the chord distances satisfy $$|x-y|
\geq |s-t|.  $$
\label{schursthm}
\end{lem}

\begin{pf} See Schur's Theorem in \cite{Ch}.
\qed
\end{pf}

We also will use results on thickness of knots developed
in
\cite{LSDR}.  Let $K$ be a smooth knot in $\Rset^3$, a
simple closed curve that is at least $C^2$ smooth.  We
assume $C^2$ smoothness throughout this paper, and the
work in  \cite{LSDR} assumed that as well; but in fact
the definition and results there can be modified to deal
with $C^{1,1}$ curves (this has been noted by R.~Litherland,
O.~Durumeric, and \cite{maddocks,cksminimum}).

For each $x \in K,$ let
$D(x,r)$ denote the  disk of radius $r$ centered at $x$ and
orthogonal at $x$ to $K$.  For sufficiently small $r>0$,
the disks are pairwise disjoint, their union forming a
tubular neighborhood of
$K$.  We define the {\it injectivity radius} of $K$,
$R(K)$, to be the supremum of such good radii. The radius
$R(K)$ measures the maximum thickness of ``rope'' that
could be used to form the curve
$K$.  Of course, $R(K)$ changes with scale.  We define the
scale-invariant {\it ropelength} or {\it length energy}
of $K$ to be
 
$$E_L(K) =\frac {\mathrm{arclength}(K)}{R(K)}\;\;.$$

The radius $R(K)$ is affected by curvature and by points
of $K$ that are far apart in the sense of arclength but
close in space.  For the latter kinds of points,
distances will be minimized at pairs of points $(x,y)$
that are critical points for the distance function
$|x-y|$. Specifically, define the  {\it critical
self-distance} of
$K$  to be the minimum of $|x-y|$ over all pairs
$(x,y) \in K \times K$, $x\not = y$, for which the chord
$x-y$ is perpendicular to $K$ at either or both of its
endpoints.  Let ${\mathrm{MinRad}}(K)$ denote the minimum
radius of curvature of
$K$.  We then have from \cite{LSDR}:

\begin{lem}  The thickness of a smooth knot is bounded by
the minimum radius of curvature and  half the critical
self-distance.  In fact,
$$R(K)= \min \left\{ {\mathrm{MinRad}}(K)\;,\; \frac12\, 
\text{critical self distance}(K)\right\}.$$
\label{lsdr}
\end{lem}

The next lemma is a consequence of Lemmas
\ref{schursthm} and \ref{lsdr}.
 
\begin{lem} Suppose $K$ is a smooth knot of thickness
$R(K)=r$.  For any $x,y \in K$, such that
$\mathrm{arc}(x,y) \geq \pi r$, we must have $|x-y| \geq
2 r$.
\label{distancebound}
\end{lem}

\begin{pf} Fix $x$ and consider first the two points $y
\in K$ for which $\mathrm{arc}(x,y) = \pi r$. Since, by
Lemma
\ref{lsdr}, the curvature of $K$ is everywhere 
$\leq 1/r$, and
$\mathrm{arc}(x,y) = \pi r$, so in particular it is $\leq
\pi r$, we can apply Schur's Theorem to the arc from $x$
to such
$y$ to conclude that $|x-y|$ is at least as large as for
the corresponding points on a circle of radius $r$,
i.e.{} $|x-y| \geq 2 r$. 

Now consider the arc $Y$ of $K$ consisting of those
points  $y$ with
$\mathrm{arc}(x,y) \geq \pi r$. Let $y_0$ be a point of
$Y$ that minimizes distance to $x$.  If $y_0$ is closer to
$x$ than $2 r$ then, from the preceding paragraph, it
cannot be an endpoint of $Y$; thus it would have to  be a
critical point for the function
$|x-y|$.  But by Lemma \ref{lsdr}, any such  critical
pair $(x,y)$ has distance
$\geq 2 r$.
\qed
\end{pf}

In the next two lemmas, we obtain the linear bound for the
energy contribution of proximal pairs, that is, points for 
which $\mathrm{arc}(x,y)\leq \pi R(K)$.

\begin{lem}  For a fixed point $s$ on a circle $C$ of
radius $R$,
$$\int\nolimits_{t\in C} \frac{1}{|s-t|^2} -
\frac{1}{\mathrm{arc}(s,t)^2} = \frac{2}{\pi R}\;.$$
\label{circle}
\end{lem}

\begin{pf}  Fix $s$ on a circle $C$ of radius $R$.  Note
that if 
$\mathrm{arc}(s,t)=\theta$, then
$|s-t|^2=R^2(2- 2 \cos\theta)$. Thus,
\begin{align*}
\int\nolimits_{t\in C} \frac{1}{|s-t|^2} -
\frac{1}{\mathrm{arc}(s,t)^2} &= 
2\int\nolimits_0^{\pi
R}\frac{1}{R^2(2-2\cos(t/R))}-\frac{1}{t^2}\,\d t \\ 
&= \frac{2}{\pi R}\;.
\end{align*}
\qed
\end{pf}

\begin{lem}  If $K$ is a smooth knot in $\Rset^3$, then
$$\int\nolimits_{x\in
K}\int\nolimits_{\mathrm{arc}(x,y)\leq \pi R(K)} 
\frac{1}{|x-y|^2} - \frac{1}{\mathrm{arc}(x,y)^2} 
\leq \frac{2}{\pi}\,E_L(K)\;.$$
\label{prox}
\end{lem}

\begin{pf} We begin by rescaling $K$ to have $R(K)=1$;
note this leaves each side of the inequality unchanged. 
Then $E_L(K)$  is just the new total arclength of $K$,
which we abbreviate $L$. By Lemma \ref{lsdr}, the
curvature of $K$ everywhere is 
$\leq \frac{1}{R(K)}=1$.  For points
$x,y$ on $K$ with $\mathrm{arc}(x,y) \leq \pi$, let
$s,t$ be points on the circle, $C$, of radius $1=R(K)$,
for which
$\mathrm{arc}(s,t) = \mathrm{arc}(x,y)$.  By Schur's
Theorem (Lemma
\ref{schursthm}), we have $|x-y| \geq |s-t|$, so
\begin{equation}
\frac{1}{|x-y|^2}-\frac{1}{\mathrm{arc}(x,y)^2} \leq 
\frac{1}{|s-t|^2}-\frac{1}{\mathrm{arc}(x,y)^2} = 
\frac{1}{|s-t|^2}-\frac{1}{\mathrm{arc}(s,t)^2}\;.
\label{arccircles}
\end{equation} 
For a fixed $x$ on $K$ and  a fixed $s$ on
$C$, (\ref{arccircles}) and Lemma \ref{circle} give us
\begin{align*}
\int\nolimits_{\mathrm{arc}(x,y)\leq \pi} 
\frac{1}{|x-y|^2} - \frac{1}{\mathrm{arc}(x,y)^2} &\leq
\int\nolimits_{\mathrm{arc}(s,t)\leq \pi}
\frac{1}{|s-t|^2}-\frac{1}{\mathrm{arc}(s,t)^2} \\ 
&= \frac{2}{\pi}\,.
\end{align*} 
Thus, 
$$
\int\nolimits_{x\in
K}\int\nolimits_{\mathrm{arc}(x,y)\leq \pi} 
\frac{1}{|x-y|^2} - \frac{1}{\mathrm{arc}(x,y)^2} \leq
\int\nolimits_{x\in K} \frac{2}{\pi} = \frac{2}{\pi}\,L\;.
$$
\qed
\end{pf} 

Finally, we need a lemma about measurable subsets
of
$K$.   When we talk about ``measure'',
we mean the one-dimensional Lebesgue measure on $K$ that is
just arclength when applied to intervals.  We actually
use only Borel sets.

We want to define subsets of $K$ of specified measure that are
closest, in terms of spatial distance, to a given point
$x\in K$.  First we need a  preliminary lemma about decomposing
measurable sets.

\begin{lem} 
\label{DecomposeASet}
Suppose $W$ is a measurable subset of $K$ and
$m_1,m_2 > 0$ are numbers with $m_1+m_2 =
\mu(W)$.  Then we can partition $W$ into measurable
sets $W_1,W_2$, with $\mu(W_j)=m_j\;, j=1,2$. 
\end{lem}

\begin{pf} 
Orient $K$, fix a point $z_0 \in K$, and
consider the intersections of $W$ with intervals
 $I_z =[z_0,z]$ as $z$ ranges over $K$.  
Since the measure of a subset of an 
interval is bounded by the length of the interval,
the function
$\mu(W\cap I_z)$ is continuous in $z$. 
The
measure $\mu(W\cap I_z)$ is arbitrarily small when
$z \approx z_0$, is monotonically
nondecreasing as $z$ moves around $K$, and
eventually exceeds $m_1$ as $z$ approaches $z_0$
from the ``back''.  
Thus there must be a value of $z$
for which $\mu(W\cap I_z) = m_1$. 
\qed
\end{pf}

In the next lemma, and the proof of the main theorem, we
use the following notation.
For $x\in K$ fixed, $a \leq b$, let $S(a,b], S[a,b], S[a]
$ denote spherical shells of radius from $a$ to $b$, and
the sphere of radius $a$, all centered at $x$. 

\begin{lem}
\label{ConstructWj}
Suppose $K$ is a smooth knot with $R(K)=1$.
Fix $x \in K$ and let $Y$ be the subset of $K$ consisting
of all points $y$ with $\mathrm{arc}(x,y) \geq \pi$. 
Let $\nu_1 ,\nu_2, \dots ,\nu_m > 0$ be numbers such
that {$\nu_1 +\nu_2 + \dots +\nu_m =
L-2\pi$} (i.e.{} $\mu (Y))$.  Then there exists a sequence
of radii
$$2 \leq \rho_1 \leq \rho_2 \leq \dots \leq \rho_m$$
and a partition of $Y$ into measurable sets 
$W_1 , W_2 ,\dots ,W_m$
such that  

(a) for each $j$, $\mu(W_j) = \nu_j$, and

(b) $W_1\subseteq Y\cap S[2,\rho_1]$ and
for each $j>1$, $W_j \subseteq Y\cap
S[\rho_{j-1},\rho_j]$.
\end{lem}

\begin{pf} 
The proof is nearly straightforward.  There is a slight complication
because the knot $K$ might intersect some spheres $S[\rho]$ in sets
of positive length.

Since $R(K)=1$ and $y\in Y$ has $\mathrm{arc}(x,y)\geq \pi$, by
Lemma \ref{distancebound},
$|x-y| \geq 2$, so $Y\subseteq S[2,L/2]$.
We define the radius $\rho_m\leq L/2$ to be the smallest radius
$\rho$ for which $Y$ is contained in $S[2,\rho]$.  

We construct radii $\rho_j$ and sets $W_j$ inductively, with a slight
variation at the final step.  For the degenerate case $m=1$, of course,
$W_1=Y$ and $\rho_1=\rho_m$ is as above.

We now begin the induction.  Suppose first that 
$\mu(Y\cap S[2]) \geq \nu_1$.  By Lemma \ref{DecomposeASet}, 
we can partition $Y\cap S[2]$
into sets $W_1\cup V$, where $\mu(W_1)=\nu_1$; let $\rho_1=2$.

Now suppose $\mu(Y\cap S[2])< \nu_1$.  Let
$$R_{\text{low}}=\{\rho\,:\,\mu(Y\cap S[2,\rho])<\nu_1\}$$
and $r_1=\sup(R_{\text{low}})$.
Similarly, let
$$R_{\text{high}}=\{\rho\,:\,\mu(Y\cap S[2,\rho])\geq \nu_1\}$$
and let $r_1'=\inf(R_{\text{high}})$.
Here $R_{\text{low}}\not = \emptyset$ because $\mu(Y\cap S[2])< \nu_1$,
so $2\in R_{\text{low}}$, and
$R_{\text{high}}\not = \emptyset$ because $\rho_m\in R_{\text{high}}$.

We claim $r_1=r_1'$.  Each element of $R_{\text{low}}$ is $<$ each element
of $R_{\text{high}}$.  Thus, $r_1\leq r_1'$.  If the inequality is strict, then
for each $\rho\in(r_1,r_1')$, $\mu(Y\cap S[2,\rho])\geq \nu_1$,
which contradicts $\rho<r_1'$.

Now define $\rho_1=r_1=r_1'$.  The set $Y\cap S[2,r_1)$ is a 
monotone union of sets, each with measure $<\nu_1$.
Thus, $\mu(Y\cap S[2,\rho_1))\leq \nu_1$.  On the other hand,
$Y\cap S(r_1',\infty)$ is a monotone union of sets, each with measure
$< \mu(Y)-\nu_1$.  Thus, $\mu(Y\cap S(\rho_1,\infty))\leq
\mu(Y)-\nu_1$, so $\mu(Y\cap S[2,\rho_1])\geq \nu_1$.

Using Lemma \ref{DecomposeASet},
we extract a subset of $Y\cap S[\rho_1]$ whose
measure is whatever $\mu(Y\cap S[2,\rho_1))$ may lack to make
$\nu_1$, and let $W_1$ be the union of this set with $Y\cap S[2,\rho_1)$.

To continue the induction, let $Y_1=Y-W_1$. Note
$Y_1\subseteq Y\cap S[\rho_1,\infty)$.

Suppose $\rho_1, \cdots, \rho_{j-1}$ and $W_1, \cdots, W_{j-1}$ have
been chosen as required in the statement of the lemma and the set
$Y_{j-1}=Y-(W_1\cup \cdots \cup W_{j-1})\subseteq S[\rho_{j-1},\infty)$.
Note that since $W_1,\cdots, W_{j-1}\subseteq Y\cap S[2,\rho_{j-1}]$,
$Y\cap S(\rho_{j-1},\infty)\subseteq Y_{j-1}$.  So for any 
$\rho,\rho'$ with
$\rho_{j-1}<\rho<\rho'$, we have
$Y_{j-1}\cap S[\rho,\rho']=Y\cap S[\rho,\rho']$.

Now proceed much like the initial case.
If $j=m$, we finish by taking $W_j=Y_{j-1}$ and $\rho_m$ as above.

There are two cases: If $\mu(Y_{j-1}\cap S[\rho_{j-1}])\geq \nu_j$, then
let $W_j$ be a subset of $Y_{j-1}\cap S[\rho_{j-1}]$ of measure
$\nu_j$, $\rho_j=\rho_{j-1}$, and $Y_j=Y_{j-1}-W_j$.

Suppose now $\mu(Y_{j-1}\cap S[\rho_{j-1}])<\nu_j$.  Define
$$R_{\text{low}}=\{\rho\,:\,\mu(Y_{j-1}\cap S[\rho_{j-1},\rho])<\nu_j\}$$
and $r_j=\sup(R_{\text{low}})$.
Similarly, let
$$R_{\text{high}}=\{\rho\,:\,\mu(Y_{j-1}\cap S[\rho_{j-1},\rho])\geq \nu_j\}$$
and let $r_j'=\inf(R_{\text{high}})$.
Here $R_{\text{low}}\not = \emptyset$ because $\rho_{j-1}\in R_{\text{low}}$ and
$R_{\text{high}}\not = \emptyset$ because $\rho_m\in R_{\text{high}}$.

As in the initial case, we have $r_j=r_j'$, and we define
$\rho_j=r_j=r_j'$.  Also, 
$\mu(Y_{j-1}\cap S[\rho_{j-1},\rho_j))\leq \nu_j$ and 
$\mu(Y_{j-1}\cap S[\rho_{j-1},\rho_j]) \geq \nu_j$.

Let $W_j$ be the union of $Y_{j-1}\cap S[\rho_{j-1},\rho_j)$ with
a subset of $Y_{j-1}\cap S[\rho_j]$ of however much additional
measure is needed to reach $\nu_j$.
Finally, let $Y_j=Y_{j-1}-W_j$.  Note then $Y_j\subseteq
Y\cap S[\rho_j,\infty)$.
\qed
\end{pf}

\section{Thick knots have bounded energy}

In this section, we prove that the M\"obius energy of a
knot is bounded by the ropelength.  Our goal is
partly the theorem itself and partly the paradigm:  any
``energy'' defined in terms of inverse-square distances
should have an analogous bound, with the proof following
this model.

\begin{thm} If $K$ is a smooth knot in $\Rset^3$ then  
\begin{equation}
\E_{O4}(K) < 4.57\, \E_L(K)^{4/3}.
\label{powbound}
\end{equation}
\label{main}
\end{thm}

\begin{thm} If $K$ is a smooth knot in $\Rset^3$ then 
\begin{equation}
\E_{O4}(K) \leq \frac{1}{4}\E_L(K)^2 \;.
\label{quadbound}
\end{equation}
\label{quadboundthm}
\end{thm}

{\bf Remark}\,
For short knots, the quadratic bound is
better than the  $4/3$ power bound.  
Comparing \eqref{powbound} to \eqref{quadbound}, we see
that one needs ropelengths over $79$ before the advantage of the
lower exponent is evident.  If one uses the actual bound
\eqref{finalbound} we obtain in the proof, which is a complicated
expression dominated by a $4/3$ power term, then that bound
is lower than \eqref{quadbound} for ropelengths $> 41$.
Computer simulations suggest \cite{stasiakideal,spatial}
the only knots that can be realized
with a ropelength $\leq 41$ are the unknot and the trefoil.

\begin{pf}[of Theorems \ref{main} and
\ref{quadboundthm}] We will obtain the quadratic bound en
route to the $4/3$ power bound.

We follow an overall plan similar to \cite{BS2,BS3}, but introduce
a limit process in the main argument.
Since
$\E_{O4}$ and
$\E_L$ are both invariant under change of scale, we begin
by rescaling $K$ to have thickness $R(K)=1$, so that
$\E_L(K)$ is just the total arclength of $K$, which we
abbreviate
$L$.

The energy
$\E_{O4}(K)$ is defined as a double-integral over $K
\times K$.  We bound separately the integral over the
portion of $K  \times K$ consisting of pairs $(x,y)$ with
$\mathrm{arc}(x,y)
\leq \pi R(K) = \pi$, and the integral over the rest of $K
\times K$.

Let $$\E_{\mathrm{prox}} = \int\nolimits_{x \in
K}\int\nolimits_{\mathrm{arc}(x,y) \leq \pi}
\frac{1}{|x-y|^2} - \frac{1}{\mathrm{arc}(x,y)^2}  
\;\;,$$ 
$$\E_{\mathrm{dist}} = \int\nolimits_{x \in
K}\int\nolimits_{\mathrm{arc}(x,y) \geq \pi}
\frac{1}{|x-y|^2} \;\;,$$ and
$$\E_{\mathrm{reg}} = \int\nolimits_{x \in
K}\int\nolimits_{\mathrm{arc}(x,y) \geq \pi}
\frac{1}{\mathrm{arc}(x,y)^2}  
\;\;.$$

So we have 
$$\E_{O4}=\E_{\mathrm{prox}}+\E_{\mathrm{dist}}-\E_{\mathrm{reg}}\;.$$

By Lemma \ref{prox}, 
\begin{equation*}
\E_{\mathrm{prox}} \leq \frac{2}{\pi}\,L\;.
\label{eprox}
\end{equation*} Also, by symmetry of the circle, 
\begin{equation*}
\E_{\mathrm{reg}}=
\int_{x \in K} \left(
2\int_{\pi}^{\frac{L}{2}}\frac{1}{t^2}\,\d t \right) \;\d
x\; = 
\;L\;\left( 2\int_{\pi}^{\frac{L}{2}}\frac{1}{t^2}\,\d t
\right) =
\frac{2}{\pi}\,L-4\;,
\label{ereg}
\end{equation*} so
\begin{equation} \label{esummary}
\E_{04} \leq E_{\mathrm{dist}}+4\;.
\end{equation}

We now bound $\E_{\mathrm{dist}}$.  We shall  bound the
inner integral and then multiply by the length of $K$ to
bound the energy.   The inner integral, for each $x$, is 
$$I^x_{\mathrm{dist}}=\int\nolimits_{\mathrm{arc}(x,y)
\geq \pi}
\frac{1}{|x-y|^2}\;\;.$$

Here first is the quadratic bound.

By Lemma \ref{distancebound}, and our rescaling to
$R(K)=1$, we have for each point $y \in K$,  
$$\mathrm{arc}(x,y) \geq \pi \Rightarrow |x-y|
\geq 2\;.$$ Thus,   
$$I_{\mathrm{dist}}^x\leq
\frac{1}{2^2}(L-2\pi)\;.$$  Multiplying by $L$, we get 
\begin{equation}
\E_{\mathrm{dist}} \leq
\frac{1}{4}\,L^2-\frac{\pi}{2}\,L\;.
\label{edistquad}
\end{equation} Combine (\ref{edistquad}) with
(\ref{esummary}) to complete a quadratic polynomial bound:
\begin{equation}
\label{quadraticbound}
\E_{O4}(K)
\leq \frac{1}{4}  L^2 - \frac{\pi}{2} L + 4\;.
\end{equation}
Because $K$ is a smooth closed curve with maximum
curvature $\leq 1$ (by Lemma \ref{lsdr} and the fact
that we have normalized $K$ to have thickness radius = 1),
the total curvature of
$K$ is at most
$L$. But by Fenchel's theorem, the total curvature of a
closed curve is at least $2\pi$. Thus $L\geq 2\pi>8/\pi$,
so the linear part of 
(\ref{quadraticbound}) is negative, and we have
\begin{equation}
\label{finalquadraticbound}
\E_{O4}(K)
\leq \frac{1}{4}  L^2\;.
\end{equation}

We now develop the  $4/3$  power bound.   
If $L<\frac{104}{3}+2\pi\approx 41$, i.e.{}
$L-2\pi<\frac{4}{3}\left(3^3-1^3\right)$, then
our proof stops here.  The quadratic bound \eqref{quadbound} is
valid and certainly $\eqref{quadbound} < \eqref{powbound}$,
so \eqref{powbound} is valid as well.
We continue under the assumption that $E_L(K)\geq \frac{104}{3}+2\pi$.

The first  observation is that for any $\rho >0$, if
$|x-y|\geq
\rho$, then the integrand $\frac{1}{|x-y|^2} \leq
\frac{1}{\rho^2}$.   In obtaining the quadratic bound,  we
stopped here, allowing the idea that with respect to each
point
$x \in K$, the whole knot (except for the arc around $x$
of length
$\pi$ in each direction) lies just at distance $2$ from 
$x$.  But {\it the knot is thick}\,:  a piece
$w$ of $K$ of length $\ell(w)$ carries along with it a
solid tube
$W$ of volume $\pi
\ell(w)$ (we still are assuming
$R(K)=1$).  Furthermore,  such a tube $W$ cannot intersect
any other part of $K$ nor any of the rest of the tube
around $K$.  This restricts how much length of $K$ can be
packed within any given distance from $x$.

Fix $x$ on $K$.   
Let
$Y_{\mathrm{dist}}$ denote the set on which we integrate
to compute
$I_{\mathrm{dist}}^x$,
that is $Y_{\mathrm{dist}}=\{y
\in K \,:\,\mathrm{arc}(x,y) \geq \pi\} $.  The arclength of
$Y_{\mathrm{dist}}$ is just $L-2 \pi$. Let
$S(\rho,\sigma]$ denote the half-open spherical shell
centered at
$x$ with radius from $\rho$ to $\sigma$. 
  Let 
$Y(\rho, \sigma]$ = $Y_{\mathrm{dist}}
\cap S(\rho, \sigma]$,  a measurable set, and let
$\ell (\rho,\sigma]$ denote its 1-dimensional measure,
i.e.{} the  measure of the length of $K$ lying in 
$S(\rho,\sigma]$.  Let 
$T(\rho,\sigma]$ be the solid tube with axis
$Y(\rho,\sigma]$, that is, the portion of the unit radius
tube about $K$ consisting of disks centered at points
$y\in Y(\rho,\sigma]$.  The 3-dimensional measure, i.e.{}
volume, of $T(\rho,\sigma]$ is
$\pi \ell(\rho,\sigma]$.
Define $S[\rho,\sigma]$, $Y[\rho,\sigma]$, $\ell[\rho,\sigma]$, and
$T[\rho,\sigma]$ analogously.
Note that $Y_{\mathrm{dist}} \subset S[2,L/2]  $. 
The subset of $Y_{\mathrm{dist}}$ 
lying at distance exactly $2$ from $x$ plays 
a special role in the analysis, so we also define  
$S[2]$, $Y[2]$, $\ell[2]$, and $T[2]$.

Let $P=(\frac{3}{4}L-\frac{3\pi}{2}+1)^{1/3}-1$.  
The constant $P$ is the
radius required so that a spherical shell from 
$r=1$ to $r=P+1$ has the same volume
as a tube of radius $1$ about a set of length $L-2\pi$.

We bound the total amount of length that can lie in each spherical
shell by a function $\ell^*$ (see \eqref{llessstar} below).
Let $\ell^*[2]=\frac{4}{3}(3^3-1^3)=\frac{104}{3}$.  
For $2\leq a<b$, let
$$\ell^*(a,b]=\frac{4}{3}\left( (\min\{b,P\}+1)^3 - (\min\{a,P\}+1)^3 
\right).$$  
In particular, when
$2\leq a<b\leq P$, this simplifies to 
$$\ell^*(a,b]=\frac{4}{3}\left( (b+1)^3 - (a+1)^3 \right).$$  
Let 
$$\ell^*[2,b] = \ell^*[2] + \ell^*(2,b].$$
For $b\leq P$, $\ell^*[2,b]=\frac{4}{3}\left((b+1)^3-1^3\right)$.
For $b>P$, notice that $\ell^*(2,b]=
\ell^*(2,P]$ so $\ell^*[2,b]=\ell^*[2,P]=L-2\pi$.

We now show that 
\begin{equation}
\label{llessstar}
\ell[2,r]\leq \ell^*[2,r] \text{ for all }2\leq r\;.
\end{equation}
The tube $T[2,r]$ has volume
$\pi\ell[2,r]$.  
But $T[2,r]$ lies entirely within $S[1,r+1]$,
whose volume is $\frac{4}{3}\pi\left( (r+1)^3-1^3\right)$.
Thus, when $2\leq r\leq P$,
\begin{align*}
\pi\ell[2,r]&\leq \frac{4}{3}\pi\left((r+1)^3-1^3\right)\\
&= \pi\ell^*[2,r].
\end{align*} 
If $r>P$, then $\ell^*[2,r]=\ell^*[2,P]=
L-2\pi\geq \ell[2,r]$. Thus, we have \eqref{llessstar}.

We next partition $Y_{\text{dist}}$ into sets $W_j$ whose distance
to $x$ is controlled, and express $I^x_{\text{dist}}$ as the sum
of the contributions from these sets.

Let $n\in\mathbb{N}$ and $\delta=\frac{P-2}{n}$.  Then,
\begin{align*}
\ell^*[2,P]&=\ell^*[2] + \ell^*(2,2+\delta] + 
\ell^*(2+\delta,2+2\delta] + \cdots + 
\ell^*(P-\delta,P]\\
&=L-2\pi.
\end{align*}

We define a sequence of sets $\{W_j\}$ using Lemma \ref{ConstructWj}
as follows.
Let 
$\nu_0=\ell^*[2]$ and 
$\nu_j=\ell^*(2+(j-1)\delta,2+j\delta]$ for $1\leq j\leq n$.  Note
that $\nu_n=\ell^*(2+(n-1)\delta,2+n\delta]=\ell^*(P-\delta,P]$.
Since $\nu_0+\nu_1+\cdots+\nu_n=L-2\pi=\mu(Y_{\mathrm{dist}})$,
by Lemma \ref{ConstructWj}, there exists a sequence
$2\leq\rho_0\leq\rho_1\leq\cdots\leq\rho_n$ and a partition
of $Y_\mathrm{dist}$ into measurable sets $W_0,\cdots,W_n$ such that:
for $1\leq j\leq n$, we have
$\mu(W_j)=\nu_j$ and $W_j\subseteq Y_\mathrm{dist}\cap S[\rho_{j-1},\rho_j]$;
for $j=0$, we have $\mu(W_0)=\nu_0$ and 
$W_0\subseteq Y_\mathrm{dist}\cap S[2,\rho_0]$.

We want to know how far each $W_j$ is from $x$, so we next bound the radii
$\rho_j$ (from below) in terms of $j$ and $\delta$.
If $y\in W_0$, then $|x-y|\geq 2$ by Lemma \ref{distancebound}.
For $1\leq j\leq n$, we claim $|x-y|>2+(j-1)\delta$ for all
$y\in W_j$, except perhaps a set of measure $0$.
Suppose $M$ is a non-zero measure subset of $W_1$ such that for all
$y\in M$, we have $|x-y|\leq 2$. Then 
$W_0\cup M\subseteq Y[2]$ and 
$\ell[2]\geq \mu(W_0)+\mu(M)> \mu(W_0)=\ell^*[2]$, which 
contradicts \eqref{llessstar}.  We continue using induction.  
Suppose for all
$0\leq k< j\leq n$, we have $|x-y|> 2+(k-1)\delta$.  Suppose
$M$ is a non-zero measure subset of $W_j$ such that for all 
$y\in M$, we have $|x-y|\leq 2+(j-1)\delta$.  Then
$W_0\cup W_1\cup \cdots \cup W_{j-1} \cup M\subseteq Y[2,2+(j-1)\delta]$.
Thus,
\begin{align*}
\pi \ell[2,2+(j-1)\delta] &\geq 
\mu(W_0)+\mu(W_1)+\cdots +\mu(W_{j-1})+\mu(M)\\
&> \mu(W_0)+\mu(W_1)+\cdots +\mu(W_{j-1})\\
&= \pi\left(\,
\ell^*[2]+\ell^*(2,2+\delta]+\cdots+\ell^*(2+(j-2)\delta,2+(j-1)\delta]
\,\right)\\
&= \pi\ell^*[2,2+(j-1)\delta],
\end{align*}
which contradicts \eqref{llessstar}.

Since $Y_{\mathrm{dist}}=W_0\cup W_1\cup \cdots \cup W_n$,
\begin{align}
I^x_{\mathrm{dist}} &= \int_{y\in W_0}\frac{1}{|x-y|^2} + 
\int_{y\in W_1}\frac{1}{|x-y|^2} + \cdots +
\int_{y\in W_n}\frac{1}{|x-y|^2}\nonumber \\
&< \frac{1}{2^2}\ell^*[2] + \frac{1}{2^2}\ell^*(2,2+\delta]+\cdots
+\frac{1}{(P-\delta)^2}\ell^*(P-\delta,P]\label{ixdistbound}.
\end{align}
The first term of \eqref{ixdistbound} is $\frac{104}{12}$.
The rest is
\begin{align*}
& \sum_{j=1}^n \frac{1}{(2+(j-1)\delta)^2}\ell^*(2+(j-1)\delta,2+j\delta)
\text{\;, which equals}\\
& \sum_{j=1}^n \frac{1}{r_j^2}\left(4(r_j+1)^2\Delta r +
4(r_j+1)\Delta r^2 + \frac{4}{3}\Delta r^3\right),
\end{align*}
where $r_j=2+(j-1)\delta$ and $\Delta r=\delta$.  This
is a Riemann sum for 
$\int_2^P\frac{4(r+1)^2}{r^2}\d r$ plus terms of higher order in $\Delta r$
whose contribution approaches $0$ as $\Delta r\to 0$.  The bound
\eqref{ixdistbound} holds for all choices of $n$, so as
$n\to \infty$, we get  
\begin{align*}
I^x_{\mathrm{dist}} 
&\leq \frac{104}{12}+\int_2^P \frac{4(r+1)^2}{r^2}\d r\\
&= \frac{8}{3}+4P-\frac{4}{P}+8\ln(P)-8\ln(2)\;.
\end{align*}
We then multiply by $L$ to get
\begin{equation}
\label{edistbound}
E_{\mathrm{dist}}\leq L
\left(\frac{8}{3}+4P-\frac{4}{P}+8\ln(P)-8\ln(2)\right).
\end{equation}
Combining \eqref{esummary} with \eqref{edistbound} yields
\begin{equation}
\label{finalbound}
E_{O4}(K)\leq L
\left(\frac{8}{3}+4P-\frac{4}{P}+8\ln(P)-8\ln(2)\right)+4\;.
\end{equation}
We substitute $P=\left(\frac{3}{4}L-\frac{3\pi}{2}+1\right)^{1/3}-1$
to obtain our final bound.

This bound is messy, but on the order of $L^{4/3}$.  We use
a computer algebra system to plot the bound divided by $L^{4/3}$.
The ratio achieves its maximum of $\approx 4.5626$ near $L=1115$.
Thus, for all $L$ we have 
$$E_{O4}(K)\leq 4.57L^{4/3}.$$
\qed
\end{pf}

For some values of $L$, in particular as
$L$ gets very large, the coefficient $4.57$ can be further reduced.
If $\frac{104}{3}+2\pi< L < 128$ or
$L> 376,\!000$,
we have $E_{O4}(K)\leq 4L^{4/3}$.  As $L$ tends to infinity,
the constant decreases to $3^{\frac{1}{3}}4^{\frac{2}{3}}\approx 3.63$.

\bibliographystyle{elsart-num}
\bibliography{newBib}

\end{document}